\documentclass[11pt]{article}
\usepackage{graphicx}
\usepackage{amscd, amssymb}
\usepackage{enumerate}
\usepackage{hyperref}
\newenvironment{eq}{\begin{equation}}{\end{equation}}
\newtheorem{proposition}{Proposition}

\newtheorem{definition}[proposition]{Definition}
\newtheorem{theorem}[proposition]{Theorem}

\newtheorem{pretheorem}[proposition]{Pretheorem}

\newcommand{\llabel}[1]{\label{#1}}
\newcommand{\comment}[1]{}
\newcommand{\sr}{\rightarrow}

\newcommand{\uu}{\underline}

\newcommand{\wt}{\widetilde}

{\vskip
3mm}

\begin{document}
\parskip = 2mm
\begin{center}
{\bf\Large The equivalence axiom and univalent models of type theory.}\\
\vskip 3mm
{\bf (Talk at CMU on February 4, 2010)}
\vskip 5mm
{\bf By Vladimir Voevodsky}

\bibliographystyle{alpha}

\begin{abstract} I will show how to define, in any type system with dependent sums, products and Martin-Lof identity types, the notion of a homotopy equivalence between two types and how to formulate the Equivalence Axiom which provides a natural way to assert that "two homotopy equivalent types are equal". I will then sketch a construction of a model of one of the standard Martin-Lof type theories which satisfies the equivalence axiom and the excluded middle thus proving that M.L. type theory with excluded middle and equivalence axiom is at least as consistent as ZFC theory. 

Models which satisfy the equivalence axiom are called univalent. This is a totally new class of models and I will argue that the semantics which they provide leads to the first satisfactory approach to type-theoretic formalization of mathematics. 
\end{abstract}

\end{center}

\subsection{Formal deduction systems and quasi-equational theories. Type systems and conservative extensions of the theory of contextual  categories.}

I will speak about type systems. It is difficult for a mathematician since a type system is not a mathematical notion. I will spend a little time explaining how I see "type systems" mathematically.  

"Type systems" are formal deduction systems of particular "flavor". So let me start with the following:

{\bf Thesis 0.} Any formal deduction system can be specified in the form of a  quasi-equational theory.

Quasi-equational theories are multi-sorted algebraic theories whose operations are given together with an ordering and "domains of definitions" which are specified by equations involving preceding operations. It is a very important class of theories.  In particular, algebraic theories (e.g. groups with given relations) are quasi-equational. The classic example of a properly quasi-equational theory is the theory of set-level categories i.e. categories up to an isomorphism.  A good source of information related to such theories is \cite{Palmgren1}.

The most important fact about these theories is the following:

{\bf Fact 1}. Any quasi-equational theory has an initial model.

With respect to the correspondence of Thesis 0, that which we call expressions, theorems, contexts etc. of a formal deduction system are elements of different sorts of the initial model of the corresponding quasi-equational theory.  We will see how this works for type systems. 

This view of formal deduction systems has many advantages. One is that it suggests a uniform approach to the formal description of various deductive systems. Another one is that "interpretations" of the deductive system are directly connected with the models of the corresponding quasi-equational theory. 

To explain how this works in the case of type systems I need to  start with the notion of a "contextual category" which was introduced by John Cartmell in \cite{Cartmell1} and studied in detail by Thomas Streicher in \cite{Streicher}:

 A contextual category is a pair of sets $C_0$, $C_1$ together with the following additional structure:
\begin{enumerate}
\item a category structure with $C_0$ being the set of objects and $C_1$ the set of morphisms (arrows),
\item an element $pt\in C_0$,
\item a map $ft:C_0\sr C_0$,
\item for each $X\in C_0$ a morphism $p_X:X\sr ft(X)$,
\item for each $X\in C_0$ and $f:Y\sr ft(X)$ an object $f^*X$ and a morphism $q(f,X):f^*X\sr X$ such that $ft(f^*X)=Y$ and the square
$$
\begin{CD}
f^*X @>q(f,X)>> X\\
@Vp_{f^*X}VV @VVp_XV\\
Y @>f>> X
\end{CD}
$$
is a pull-back square.
\end{enumerate}
These data should satisfy the following conditions:
\begin{enumerate}
\item $ft(pt)=pt$,
\item $Id^*_{ft(X)}=X$ and $q(Id_{ft(X)},X)=Id_X$,
\item for $Z\stackrel{g}{\sr}Y\stackrel{f}{\sr}ft(X)$ one has $(fg)^*(X)=g^*f^*(X)$ and $q(fg,X)=q(f,X)q(g,f^*X)$.
\end{enumerate}

Note that contextual categories are "set-level" entities i.e. we consider two contextual categories to be "the same" if they are isomorphic to each other. We do not introduce any notion of an equivalence for contextual categories. 

Defined as above contextual categories are models of a finitary quasi-equational theory $\cal CC$. It is shown in \cite{Palmgren1} that axioms of the form "there exists a unique ..." can be rewritten in quasi-equational form which allows us to include the condition the canonical squares in contextual categories are pull-back squares.  

Both Cartmell and Streicher add to the list of axioms for contextual categories the axioms $(ft(X)=X)\Rightarrow (X=pt)$ and $\forall X\in C_0,\exists n\ge 0, ft^n(X)=pt$. I do not include these because the first one is never used and the second one is not a quasi-equational one. 

Let me connect now type systems with quasi-equational theories and contextual categories. 

{\bf Thesis 1}. Any type system defines a contextual category whose objects are equivalence classes, with respect to the definitional equality, of valid contexts of  the type system  and morphisms from a context $(y_1:{\bf S}_1,\dots,y_m:{\bf S}_m)$ to a context $(z_1:{\bf T}_1,\dots,z_n:{\bf T}_n)$ are (equivalence classes of) sequences of term sequents of the form 
$$(\vec{y}:\vec{\bf S}\vdash {\bf u}_1:{\bf T}_1),\,\,\,(\vec{y}:\vec{\bf S}, z_1:{\bf T}_1\vdash {\bf u}_2:{\bf T}_2),\,\,\,\dots,$$
$$\dots\,\,\,(\vec{y}:\vec{\bf S},z_1:{\bf T}_1,\dots,z_{n-1}:{\bf T}_{n-1}\vdash {\bf u}_n:{\bf T}_n)$$

{\bf Thesis 2.} The contextual categories corresponding to most type systems are the initial model of quasi-equational theories which are conservative (without adding new sorts) extensions of $\cal CC$.

In view of Thesis 2, one can view type systems as "syntactic implementations" of initial models of quasi-equational theories. Therefore, if I want to construct a "model" of a type system in some other  mathematical theory (e.g. ZFC) I should solve two independent problems:

1. Write down a description of the quasi-equational theory whose initial model is supposed to be "implemented" by the type system and prove that the type system does indeed implement such a model. This involves working with a formal description of the syntax and the reduction rules of the type system.  

2. Construct, in the framework of my mathematical theory, a model of this quasi-equational theory with the desired properties.

This strategy is, as far as I know, new. I very much hope that the future work on the semantics of type system will follow it or some version of it because the current constructions relating type systems to more familiar mathematical entities are very difficult to follow. I also suggest that in the development of future type systems we start with providing a description of the quasi-equational theory whose initial model the type system is supposed to implement, followed by the description of the syntactic implementation (i.e. the type system per-ce) and by the proof that the suggested syntax indeed defines an initial model.  In the future this process, I hope, will be formalized and automated. This raises the following interesting question:

{\bf Question 1.} What is the weakest proof environment which can be used to construct formal proofs connecting formal deduction systems (in particular type systems) to initial models of quasi-equational theories?

In the rest of my talk I will mostly address the second step of this strategy i.e. the construction, in the framework of ZFC, of contextual categories with additional structures (e.g. the ones corresponding to dependent products, dependent sums, unverses, Martin-Lof equality etc.)

\subsection{Contextual categories defined by universes (morphisms) in lcccs. Structures on universes corresponding to structures on contextual categories corresponding to standard rules of type systems.}

The next idea which I want to propose concerns a new method of construction of contextual categories with different structures starting with objects "of level one" i.e. object which are defined only up to an equivalence. 

Let $\cal C$ be a category. Define a universe in $\cal C$ as a morphism $p:\wt{U}\sr U$ together with a choice, for any morphism $f:X\sr U$, of a pull-back square of the form
$$
\begin{CD}
(X;f) @>Q(f)>> \wt{U}\\
@Vp_{(X,f)}VV @VVpV\\
X @>f>> U
\end{CD}
$$
and with a choice of a final object $pt$ in $\cal C$.

We will write $(X;f_1,\dots,f_n)$ for $((X;f_1,\dots,f_{n-1});f_n)$.  Given a universe $p$ in $\cal C$ define a contextual category $CC=CC({\cal C},p)$ as follows. Objects of $CC$ are sequences of the form $(F_1,\dots,F_m)$ where $F_1:pt\sr U$ and $F_{i+1}:(pt;F_1,\dots,F_i)\sr U$. Morphisms from $(F_1,\dots,F_m)$ to $(G_1,\dots,G_n)$ are morphisms from $(pt;F_1,\dots,F_m)\sr (pt;G_1,\dots,G_n)$ in ${\cal C}$. The rest of the contextual structure is easy to infer. 

Note that up to an equivalence, $CC$ is just a full subcategory in ${\cal C}$ which consists of objects $X$ such that the morphism $X\sr pt$ is a finite composition of morphisms which can be "induced" from $p$. However, the main good property of this construction is the following:
\begin{proposition}
The contextual category $CC({\cal C},p)$ is well defined up to a canonical {\em isomorphism} by the equivalence class of $({\cal C}, p)$. 
\end{proposition}

In particular, up to a canonical isomorphism this contextual category does not depend on the choice of $pt$ or the pull-back squares which is why we denote it simply by $CC({\cal C},p)$. 

The most important, for type theory, structure on contextual categories is the "products of families of types" (cf. \cite[p. 71]{Streicher}). This structure, in type theoretic terms, corresponds to having dependent products of the form
$$
\frac{\Gamma, a:{\bf F}, b:{\bf G}\vdash}{\Gamma, c:\prod a:{\bf F}.{\bf G}\vdash}
$$
with the usual $\lambda$ and $eval$ and both the $\beta$  and the $\eta$ reductions. 

For the following proposition recall that an lccc is a locally cartesian closed category i.e. a category $\cal C$ such that all its slice categories ${\cal C}/X$ have finite products and internal hom-objects.
\begin{proposition}
Let $\cal C$ be an lccc and $p:\wt{U}\sr U$ a universe in $\cal C$. Then any choice of two morphisms 
$$\wt{P}:\uu{Hom}_{U}(\wt{U}, \uu{U}\times \wt{U})\sr \wt{U}$$
$$P:\uu{Hom}_{U}(\wt{U},\uu{U}\times U)\sr U$$
such that the square
$$
\begin{CD}
\uu{Hom}_{U}(\wt{U}, \uu{U}\times \wt{U}) @>\wt{P}>> \wt{U}\\
@V\uu{Hom}_{U}(\wt{U}, Id_U\times p)VV @VVpV\\
\uu{Hom}_{U}(\wt{U},\uu{U}\times U) @>P>> U
\end{CD}
$$
is a pull-back square, defines on $CC({\cal C},p)$ a "products of families of types" structure which is compatible in the obvious sense with the forgetting functor to ${\cal C}$.
\end{proposition}

In the proposition above and in diagrams below I underline the copy of $U$ the projection to which is used to consider the product as an object over $U$. 

A slightly more complex structure on $p$ defines the structure on $CC({\cal C},p)$ corresponding to the dependent sums in type theory. Type theoretic universes, their mutual relations $\le $ and $:$ their property of being closed under dependent products and sums etc. also correspond to relatively simple structures on $p$.  

To get on $CC({\cal C},p)$ a structure corresponding to the Martin-Lof equality in the form
$$\frac{\Gamma, t_1:{\bf T}, t_2:{\bf T}\vdash}{\Gamma, u:Eq({\bf T},t_1,t_2)\vdash}$$
one needs a pair of morphisms $\Omega:\wt{U}\sr \wt{U}$ and $Eq:\wt{U}\times_U \wt{U}\sr U$ such that the square
\begin{eq}
\llabel{2010.2.2.eq1}
\begin{CD}
\wt{U} @>\Omega>> \wt{U}\\
@V\Delta_{\wt{U}}VV @VVpV\\
\wt{U}\times_U \wt{U} @>Eq>> U
\end{CD}
\end{eq}
commutes and a section of the morphism
$$\uu{Hom}_U(E \wt{U},\uu{U}\times U)\sr \uu{Hom}_U(E \wt{U},\uu{U}\times U)\times_{\uu{Hom}_B(\wt{U},\uu{U}\times U)}\uu{Hom}_U(\wt{U},\uu{U}\times \wt{U})$$
which is defined by the commutative square
$$
\begin{CD}
\uu{Hom}_U(E \wt{U}, \uu{U}\times \wt{U}) @>\uu{Hom}_B(i,\uu{U}\times \wt{U})>> \uu{Hom}_U(\wt{U},\uu{U}\times \wt{U})\\
@V\uu{Hom}_U(E \wt{U},Id_U\times p)VV @VV\uu{Hom}_U(\wt{U},Id_U\times p)V\\
\uu{Hom}_U(E \wt{U}, \uu{U}\times U) @>\uu{Hom}_B(i,U)>> \uu{Hom}_U(\wt{U},\uu{U}\times U)
\end{CD}
$$
where $E \wt{U}=(\wt{U}\times_U \wt{U})\times_U \wt{U}$ (considered over $U$ through the projection $\wt{U}\times_U\wt{U}\sr U$) and $i:\wt{U}\sr E \wt{U}$ is the morphism defined by (\ref{2010.2.2.eq1}).  

Such a section defines for any $f:X\sr U$ and any commutative square of the form
$$
\begin{CD}
(X;f) @>>> \wt{U}\\
@VId_X\times_U iVV @VVpV\\
E X @>>> U
\end{CD}
$$
where $EX=X\times_U E\wt{U}$, a morphism $E X\sr \wt{U}$ making the two triangles commutative. 

There is a lot of very interesting questions related to what is "the most general" inductive construction in type theory and how to describe structures on universe maps $p$ which generate operations of the corresponding form on contextual categories $CC({\cal C},p)$. For now I will leave these questions open.

So now we have two ingredients - the interpretation of type systems as implementations of the initial models of conservative extensions of the theory of contextual categories and a construction which produces models of such extensions from lcccs with a morphism $p:\wt{U}\sr U$ and some structures on this morphism.

To describe the models of type theory which I call univalent models it remains to specify a particular lccc and a particular morphism $p:\wt{U}\sr U$ which has all the required structures. This is where the homotopy theory becomes essential.

\subsection{Simplicial sets. Universal Kan fibrations as a universe. Univalence property.} 

Let me recall some basic definitions. The simplicial category $\Delta$ is the category whose objects are natural numbers (denoted $[n]$) and morphisms from $[m]$ to $[n]$ are order preserving maps from the finite set $\{0,\dots,m\}$ to the finite set $\{0,\dots,n\}$. A simplicial set is a contravariant functor from $\Delta$ to $Sets$. The category of simplicial sets is denoted $\Delta^{op}Sets$. It is a Grothendieck topos and in particular an lccc. 

While the particulars of $\Delta^{op}Sets$ will be used below one can with more or less work replace it by any other sufficiently "good" model for the homotopy theory. The fact that the localization of $\Delta^{op}Sets$ with respect to the class $W$ of weak equivalences which I will define in a few moments, is equivalent to the category of "nice" topological spaces and homotopy classes of maps is undoubtedly one of the most important and most useful, discoveries of the 20-ies century mathematics.  Because of this fact one can, when working up to homotopy, think of simplicial sets as of combinatorial representations of shapes or spaces, of simplicial paths as paths on these spaces etc. This provides a lot of the underlying intuition for the univalent models. For a classic introduction including minimal fibrations see \cite{Ma}.

One defines $\Delta^n$ as the functor represented by $[n]$. Its geometrical realization is the usual $n$-simplex 
$$\Delta^n_{top}=\{x_0,\dots,x_n\,|\,x_0,\dots,x_n\ge 0,\,\,\,\sum_i x_i=1\}.$$
 Let $\Lambda^n_k$, $k=0,\dots,n$ be the simplicial subset of $\Delta^n$ which is the boundary $\partial\Delta^n$ of $\Delta^n$ without the $k$-th face. 

One defines Kan fibrations as morphisms $q:E\sr B$ such that for any $n$ and $k$ and any commutative square of the form
$$
\begin{CD}
\Lambda^n_k @>>> E\\
@VVV @VVqV\\
\Delta^n @>>> B
\end{CD}
$$
where the left had side vertical arrow is the natural inclusion, there exists a morphism $\Delta^n\sr E$ which makes both triangles commutative.  A simplicial set whose morphism to the point is a Kan fibration is called a Kan simplicial set. 

A morphism $i:A\sr X$ is called an anodyne morphism if for any Kan fibration $q:E\sr B$ and any commutative square of the form
$$
\begin{CD}
A @>>> E\\
@ViVV @VVqV\\
X @>>> B
\end{CD}
$$
there exists a morphism $X\sr E$ which makes both triangles commutative. Clearly, the inclusions $\Lambda^n_k\sr \Delta^n$ are anodyne but there are actually many more anodyne morphisms.  A morphism which becomes an isomorphism if one inverts all anodyne morphisms is called a weak equivalence. A simplicial set whose projection to the point is a weak equivalence is called (weakly) contractible. 

Let now $\alpha$ be a sufficiently large cardinal. For simplicity less us assume that it is a very large cardinal - a strongly inaccessible one in some appropriate sense. Let us say that an ordered morphism $f:X\sr Y$ between two simplicial sets is a morphism together with well-orderings of all fibers. Since there is at most one isomorphism between two well ordered sets there is at most one codomain fixing, order preserving isomorphism between two ordered morphisms. This makes it possible to consider the universal ordered morphism $\wt{V}_{<\alpha}\sr V_{<\alpha}$ with fibers of cardinality $<\alpha$ such that any  ordered morphism $f:Y\sr X$ with such fibers can be included into a unique pull-back square of the form
$$
\begin{CD}
Y @>>> \wt{V}_{\alpha}\\
@VVV @VVV\\
X @>>> V_{<\alpha}
\end{CD}
$$
Since the class of Kan fibrations is closed under pull-backs (this is obvious from their definition) there is a sub-morphism $p_{<\alpha}:\wt{U}_{<\alpha}\sr U_{<\alpha}$ in $\wt{V}_{<\alpha}\sr V_{<\alpha}$ which is the universal ordered Kan fibration with fibers of cardinality $<\alpha$. 
\begin{theorem}
\llabel{2010.2.2.th1}
For a sufficiently large $\alpha$ there are structures on $p_{<\alpha}$ corresponding to the products of families of types, dependent sums, impredicative Prop and Martin-Lof equalities. Moreover, in a well defined sense, the space of such structures on $p_{<\alpha}$ is contractible.
\end{theorem}

Note that the second half of the theorem asserts that  these structures on $p_{<\alpha}$ are essentially unique. There are also structures corresponding to universes. As with the inductive types (and more so) I am not completely comfortable with  my current understanding of these "universe structures" and will leave this issue open. 

\begin{definition}
The models of type theory which correspond to the contextual categories $CC(\Delta^{op}Sets, p_{<\alpha})$ are called standard univalent models. 
\end{definition}

These models were foreseen to some extend by Hofmann and Streicher in \cite{Hofmann1}. They have also foreseen there some of the implications which this view of the semantics will have on the formalization of mathematics.    

Two properties of $p_{<\alpha}$ are of the key importance for proving Theorem \ref{2010.2.2.th1}. Both are proved using the theory of so called minimal fibrations. The first one is simple:
\begin{theorem}
For any infinite $\alpha$ the simplicial set $U_{<\alpha}$ is a Kan simplicial set.
\end{theorem}
The second one requires an additional definition. It also brings us back to the equivalence axiom. 

For any morphism $q:E\sr B$ consider the space (simplicial set) $\uu{Hom}_{B\times B}(E\times B, B\times E)$. If $q$ is a fibration then this space contains, as a union of connected components, a subspace $Eq(E\times B, B\times E)$ which corresponds to morphisms which are weak equivalencies and the obvious morphism from $\delta:B\sr B\times B$ to $\uu{Hom}_{B\times B}(E\times B, B\times E)$ over $B\times B$ factors uniquely through a morphism
$$m_q:B\sr Eq(E\times B, B\times E)$$
\begin{definition}
A fibration $q$ is called univalent if $m_q$ is a weak equivalence.
\end{definition}

There is a whole nice theory of univalent fibrations which I have no time to go into. It has not been studied previously I guess because univalent fibrations are not stable under base change. The key for us is the following.
\begin{theorem}
For any infinite $\alpha$, the fibration $p_{<\alpha}$ is univalent.
\end{theorem}

This theorem is essentially equivalent to the property that the models of type theory based on $p_{<\alpha}$ satisfy the {\em equivalence axiom}.  Let me explain this axiom using the language of Coq.
In this explanation I will only use the dependent sums and equality whose univalent interpretation is guaranteed by Theorem \ref{2010.2.2.th1}. It is however my current understanding that the univalent interpretation can also be extended, in an essentially unique way,  to all the inductive constructions supported by the calculus of inductive constructions which Coq is based on. I have not yet checked the mutually inductive definitions and the coinductive definitions which are also a part of Coq language. This being said, I believe that the following is true:

\begin{pretheorem}
A proof of "False" in Coq under the assumption of the equivalence, contractible choice and the excluded middle axioms would imply inconsistency of ZFC.
\end{pretheorem}

Actually proving this statement would be a good application of the general strategy for the construction of type systems and their models which I have outlined at the beginning of the lecture. 

The equivalence axiom is best explained with a piece of an actual Coq code directly in CoqIde.

\def\cprime{$'$}

\end{document}